\documentclass[10pt]{article}
\usepackage{lineno,hyperref}
\usepackage{amssymb}
\usepackage{amsmath}
\usepackage{newtxtext}       %
\usepackage[varvw]{newtxmath} 
\renewcommand{\vec}[1]{\boldsymbol{#1}}
\begin{document}
\author{Alexander I. Zhmakin\footnote{Ioffe Physical  Technical Institute \& SoftImpact, Ltd., e-mail: a.zhmakin0@gmail.com}}
\title{A Compact Introduction to Fractional Calculus}

\maketitle
\tableofcontents
\pagenumbering{arabic}

\section{Fractional Derivatives}
Fractional calculus (FC) is now an efficient tool for problems in science and engineering \cite{old74,mil93,kir94,sam87,car97,pod98,nakh03,tar11}.
The term "fractional" is kept for the historical reasons --- it is a misnomer since the order can be real \cite{gor03,ben13}.

Historical survey of the development of FC starting from the letter by Gottfried Leibniz to Guillaume l'H\^opital (1695) including contributions by Joseph Liouville, Bernhard Riemann, Niels Abel, Gr\"unwald,  Aleksey Letnikov, Gerasimov, Marcel Riesz, Magnus Mittag-Leffler, Paul L\'evy, Raoul Nigmatullin, Yuri Rabotnov, Arthur Erd\'elyi and others during the XIX and XX centuries could be found in refs. \cite{ross97,hil08,kir10,kir11,kir14}. 

Lorenzo \& Hartley \cite{lor98} analysed the minimal set of criteria for the generalized calculus formulated by B. Ross:
{{\em analyticity}: if $f (z)$ is an analytic function of the complex variable $z$, the derivative $D_z^{\nu} f (z)$ is an analytic function of $z$ and $\nu$;}
{{\em backward compatibility}: the operation $D_z^{\nu} f (z)$ must produce the same result as ordinary differentiation when $\nu = n$ is a positive integer; the operation $D_z^{\nu} f (z)$ must produce the same result as ordinary $n$-fold integration when $n$ is a negative integer;  $D_z^{\nu} f (z)$ must vanish along with its $n - 1$ derivatives at  $x = c$;}
{{\em zero property}: the operation of order zero leaves the function unchanged $D_z^0 f (z) = f (x) $;}
{{\em linearity}: the fractional operators must be linear $D_z^{\nu} [a f (x) + b g (x) ] = a D_z^{\nu} f (x) + D_z^{\nu} g (x)$;}
{{\em composition (index law)}: the law of exponents for integration of arbitrary order $D_z^{\nu} D_z^{\mu} f (x) = D_z^{\nu + \mu} f (x)$.}

The fractional derivatives are based on the extension of the repeated integration and are defined either by the continuation of the fractional integral to the negative order or by the integer order derivatives of the fractional integrals \cite{hil08}.
There is no unique definition \cite{sam87,pod98,tar11} (and notation is not standardized \cite{hil08,oli14}). 

There are several kinds of definitions of the fractional derivatives (Riemann-Liouville, Caputo, Gr\"unfeld-Letnikov, Riesz, Weyl, Marchaud, Caputo-Fabrizio, Yang, Chen, He and  others) that are not equivalent \cite{li07,li11}. The initial conditions for the Caputo derivative are expressed in terms of the initial values of the integer order derivatives; for the zero initial conditions  Riemann-Liouville, Caputo and Gr\"unwald-Letnikov derivatives coincide \cite{rah10}.
 Most of the fractional derivatives are defined through the fractional integral thus derivatives inherent some non-local behaviour \cite{iyi16}.
The following relation is valid fo all types of the fractional derivatives \cite{zas08}
\begin{equation*}
	\frac{d^{\alpha + \beta}}{d t^{\alpha + \beta}} = \frac{d^{\alpha}}{d t^{\alpha}} \frac{d^{\beta}}{d t^{\beta}} = \frac{d^{\beta}}{d t^{\beta}} \frac{d^{\alpha}}{d t^{\alpha}}.
\end{equation*}

The most frequently used are the Riemann-Liouville (e.g., in the calculus), the Caputo (e.g., in the physics, the numerical computations) and Gr\"unwald-Letnikov (e.g., in the signal processing, the engineering)  fractional derivatives \cite{agu14,oli14}.

Grigoletto \& de Oliveira \cite{gri13} considered the generalization of the fundamental theorem of calculus  --- Fundamental Theorem of Fractional Calculus (FTFC) for the cases of the Riemann-Liouville, Liouville, Caputo, Weyl and Riesz derivatives.

Baleanu \&  Fernandez \cite{bal19} considered the possible classification of the fractional operators into broad classes under some restrictions and criteria. In particularity, many operators could be considered as the special cases of \cite{fer19}
\begin{equation*}
	_c^AI_x^{\alpha,\beta} \int\limits_c^x (x - t)^{ \alpha - 1} A \left( (x - t)^{\beta} \right) f(t) d t, \qquad A (z) = \sum\limits_{k =0}^{\infty} a_k z^k
\end{equation*}
where $c$ is a constant often taken as zero or $\infty$,  $\alpha$ and $\beta$ are complex parameters with positive real parts, and $
	A (z)$ is a general analytic function whose
coefficients $a_k \in  \mathcal{C}$ are permitted to depend on $\alpha$ and $\beta$, $x$ as a real variable larger than $c$.

\subsection{Riemann-Liouville Fractional Integral}
Both the Riemann-Liouville (RL) and  the Caputo fractional derivatives are based on the RL fractional integral  that for any $\alpha > 0$ is defined as \cite{mai07,del13}
\begin{equation}
	\label{lrl}
	J_{a^+}^{\alpha} f (x) = \frac{1}{\Gamma (\alpha)} \int\limits_a^x (x - t)^{\alpha - 1} f (t) d t,
\end{equation}
If $\alpha = 0$, $J_{a^+}^0 = I$, $I$ is the identity operator.
Here 
$$\Gamma (\alpha)  = \int\limits_0^{\infty} \exp (- \alpha) u^{\alpha - 1} d u$$
is the Euler Gamma function.
This integral exists if $f (t)$ is the locally integrable function and for $t \rightarrow 0$ behaves  like $O (t^{-\nu})$ with $\nu < \alpha$. To get the strict mathematical rigor it is possible to use the framework of the Lebesgue spaces of the summable functions or the Sobolev spaces of the generalized functions \cite{gor97a}.

The RL  integral is a generalization of Cauchy's formula for an n-fold integral
\begin{equation*}
	\int\limits_a^x dx_1	\int\limits_a^{x_1} dx_2	\dots \int\limits_a^{x_{n - 1}} dx_n	= \frac{1}{(n - 1)!} \int\limits_a^x (x - t)^{n - 1} dt
\end{equation*}
using the relation 
\begin{equation*}
	(n - 1)! = \prod_{k = 1}^{n - 1} k = \Gamma (n).
\end{equation*}

The equation (\ref{lrl}) is {\em left-sided} RL integral. The {\em right-sided} RL integral is written as 
\begin{equation}
	\label{rrl}
	J_{b^-}^{\alpha} f (x) = \frac{1}{\Gamma (\alpha)} \int\limits_x^b (t - x)^{\alpha - 1} f (t) d t.
\end{equation}

The RL integral is a case of the convolution integral of the Volterra type \cite{mai00} 
\begin{equation*}
	K * f (x) = \int\limits_a^b k (x - t) f (t) dt.
\end{equation*}

The RL integral has the {\em semi-group} property (also called {\em additivity law} \cite{hil08}):
\begin{equation*}
	J_{a^+}^{\alpha} J_{a^+}^{\beta} f (x) = 	J_{a^+}^{\alpha + \beta} f (x), \quad \alpha > 0, \quad \beta > 0
\end{equation*}
which implies the {\em commutative} property \cite{mai00}: $	J_{a^+}^{\beta} J_{a^+}^{\alpha} = J_{a^+}^{\alpha} J_{a^+}^{\beta}$.	

The RL fractional integral coincides with the classical definition in the case $\alpha \in \mathcal{N}$. The fractional integration improves the  smoothness  of functions \cite{die10}.

Sometimes the RL integral could be expressed via the elementary functions, e.g., 
\begin{equation*}
	J_{a^+}^{\alpha} (x - a)^{\mu} = \frac{\Gamma (\mu + 1)}{\Gamma (\alpha + \mu +1)} (x - a)^{\alpha + \mu}.
\end{equation*}

A particular case of the RL fractional integrals is the Liouville fractional integrals \cite{gri13} that is obtained by transitions $a \rightarrow - \infty$ and $b \rightarrow  \infty$ in equations (\ref{lrl}) and (\ref{rrl}) as
\begin{equation*}
	J_{+}^{\alpha} f (x) = \frac{1}{\Gamma (\alpha)} \int\limits_{- \infty}^x (x - t)^{\alpha - 1} f (t) d t,
	\quad
	J_{-}^{\alpha} f (x) = \frac{1}{\Gamma (\alpha)} \int\limits_x^{\infty} (t - x)^{\alpha - 1} f (t) d t.
\end{equation*}

\subsection{Riemann-Liouville Fractional Derivative} 
{The left and the right Riemann-Liouville fractional derivatives are defined as \cite{pod98}
	\begin{equation}
		\label{lrld}
		D_{a^+}^{\alpha} [ f (x)]= \left \{
		\begin{array}{rl}
			\displaystyle \frac{1}{\Gamma (1 - \alpha)} 
			\frac{d}{dx}
			\int\limits_a^x (x - t)^{- \alpha} f(t) dt,  & \quad \alpha \in (0,1)\\
			\displaystyle \frac{df (x)}{dt},
			&  \quad \alpha = 1
		\end{array} \right.
	\end{equation}
	and
	\begin{equation}
		\label{rrld}
		D_{b^-}^{\alpha} [ f (x)]= \left \{
		\begin{array}{rl}
			\displaystyle \frac{1}{\Gamma (1 - \alpha)} 
			\frac{d}{dx}
			\int\limits_x^b (t - x)^{- \alpha} f(t) dt,  & \quad \alpha \in (0,1)\\
			\displaystyle \frac{df (x)}{dt},
			&  \quad \alpha = 1
		\end{array} \right.
	\end{equation}
	Operator $D_{a^+}^{\alpha}$ is left-inverse meaning that $D_{a^+}^{\alpha} J_{a^+}^{\alpha} = I$, $I$ is the identity operator. Thus $D_{a^+}^{\alpha} J_{a^+}^{\alpha} f = f$ but the unconditional semigroup property of fractional differentiation in the RL sense does not hold: Diethelm \cite{die10} gives examples where $D_{a^+}^{\alpha_1} D_{a^+}^{\alpha_2} f = D_{a^+}^{\alpha_2} D_{a^+}^{\alpha_1} f \ne D_{a^+}^{\alpha1 + \alpha_2} f$ and $D_{a^+}^{\alpha_1} D_{a^+}^{\alpha_2} f \ne D_{a^+}^{\alpha_2} D_{a^+}^{\alpha_1} f = D_{a^+}^{\alpha1 + \alpha_2} f$.
	
	Atangana \& Secer \cite{atang13} presented tables of the RL derivatives of the  trigonometric and some special functions.
	The fractional RL derivative of the power function is
	\begin{equation*}
		D_{a^+}^{\alpha} t^{\nu} = \frac{\Gamma (1 + \nu)}{\Gamma (1 + \nu - \alpha)} t^{\nu - \alpha}
	\end{equation*}
	and, particular, the derivative of a constant
	$	D_{a^+}^{\alpha} 1 = {t^{- \alpha}}/{\Gamma (1 - \alpha)}$. 
	
	Since the fractional RL derivative of a constant is not zero, thus the magnitude of the fractional derivative changes with adding of the constant. 
	
	Jumarie \cite{jum} suggested a modification to remove this drawback. He started with a fractional derivative ({\em F-derivative})
	\begin{equation*}
		f^{\alpha} (x) = \lim_{h \rightarrow 0} \frac{\Delta^{\alpha} f (x)}{h^{\alpha}}
	\end{equation*}
	based on the fractional difference $\Delta^{\alpha} f (x)$ of order $\alpha, \quad \alpha \in \mathfrak{R}, \quad 0 < \alpha \le 1$.
	
	Jumarie proposed the modification of the fractional RL derivative
	\begin{equation*}
		\displaystyle \frac{1}{\Gamma (1 - \alpha)} 
		\frac{d}{dx}
		\int\limits_a^x (x - t)^{- \alpha} (f(t) - f (0))dt.
	\end{equation*}

	\subsubsection{Leibniz' formula}
	The classical Leibnitz' formula for the first-order derivative (i.e. when $n \in \mathcal{N}$) is
	\begin{equation*}
		D^n [f (x) g (x)] = \sum_{k = 0}^n \binom{n}{k} [D^k g (x) D^{n - k} f (x)]
	\end{equation*}
	where $f (x)$ and $g (x)$ are the $n$-time differentiable functions.
	
	The fractional derivatives violate the classical Leibnitz' rule \cite{tar1,tar2}. Generalization of the Lebnitz' formula was  developed by Osler \cite{osl1,osl2}.
	
	The Leibniz' formula for the differentiation of the product of the functions for the RL operators for the functions that are analytic on $(a - h, a + h)$ is written as \cite{die10}
	\begin{equation*}
		D_{a^+}^{n} [f g] (x) = \sum_{k = 0}^{\lfloor n \rfloor} \binom n k (D_{a^+}^{k} f) (x) (D_{a^+}^{n - k} g) (x) + 
		\sum_{k = [n] + 1}^{\infty} \binom n k (D_{a^+}^{k} f) (x) (J_{a^+}^{n - k} g) (x).
	\end{equation*}
	where $\lfloor \quad \rfloor$ denotes the floor function.
		Jumarie studied the Leibniz' formula for the differentiation of the product of the non-differentiable functions \cite{jum13}. 
	
	\subsubsection{Fa\'a di Bruno formula (the chain rule)}
	For the functions $f$ and $g$ with a sufficient number of the derivatives and $n \in \mathcal{N}$ \cite{die10,del13,tar16a}
	\begin{equation*}
		D^n [g (f (\cdot))] (x) = \sum (D^k g) \prod_{m = 1}^n (D^m f (x)^{b_m}
	\end{equation*}
	where the sum is over all partitions of $\{1, 2, \dots, n\}$ and for each partition $k$ is its number of the blocks and $b_j$ is the number of the blocks with exactly $j$ elements.
	
	Tarasov \cite{tar16} analysed the simplified chain rules suggested by Jumarie \cite{jum07,jum07a,jum13a} and found that these simplifications are not universally valid.

	\subsubsection{Fractional Taylor expansion}
	The fractional Taylor expansion is written as \cite{die10,qi13,odi07}
	\begin{equation*}
		f (x) = \frac{(x - a)^{n - m}}{\Gamma (n - m +1)} \lim_{z \rightarrow a +} J_a^{m - n} f (z) +
	\end{equation*}
	\begin{equation*}	
		\sum_{k = 1}^{m - 1} \frac{(x - a)^{k + n - m}}{\Gamma (k + n - m +1)} \lim_{z \rightarrow a +} D^k J_a^{m - n} f (z) + J_a^n D_a^n f (x).
	\end{equation*}
	
	\subsubsection{Symmetrised space derivative}
	Vermeersch \& Shakouri \cite{ver14} formulated the symmetrised space derivatives of the fractional order between 1 and 2 and between 0 and 1:
	\begin{itemize}
		\item {$1 < \alpha < 2$. 			
			The symmetrised space derivative of the function $g (x)$ that is integrable over the entire real axis is
			\begin{equation*}
				\frac{\partial^{\alpha} g}{\partial |x|^{\alpha}} = \frac{\partial}{\partial x} \left[w_{\alpha} \star \frac{\partial g}{\partial x} \right] 
			\end{equation*}	
			where $\star$ denotes the convolution and $w_{\alpha}$ is an unknown kernel function with the Fourier image found to be
			$W_{\alpha} = {1}/{|\xi|^{2 - \alpha}}$.
			The Fourier inversion yields
			\begin{equation*}
				w_{\alpha} = \frac{1}{2 \pi} \int\limits_{- \infty}^{\infty} \frac{exp (j \xi x) d \xi}{|\xi|^{2 - \alpha}} = \frac{|x|^{- (\alpha - 1)}}{2 \Gamma (2 - \alpha) \cos [(2 - \alpha) \frac{\pi}{2}]}.	
			\end{equation*}
			Thus
			\begin{equation*}
				\frac{\partial^{\alpha} g}{\partial |x|^{\alpha}} =	\frac{1}{2 \Gamma (2 - \alpha) \cos [(2 - \alpha) \frac{\pi}{2}]} \frac{\partial}{\partial x} \int\limits_{- \infty}^{\infty} \frac{\frac{\partial g}{\partial x} (x^{\prime}) d x^{\prime}}{|x - x^{\prime}|^{\alpha - 1}}.
			\end{equation*}
		}
		\item {$0 < \alpha < 1$}. 	The symmetrised space derivative of the function $g (x)$ is 
		
		\begin{equation*}
			\frac{\partial^{\alpha} g}{\partial |x|^{\alpha}} = w_{\alpha} \star \frac{\partial g}{\partial x}. 
		\end{equation*}	
		
		The Fourier image of the kernel function  $w_{\alpha}$ is $W_{\alpha} = {j \cdot sgn (\xi)}/{|\xi|^{1 - \alpha}}$.	
			Performing the Fourier inversion, the authors get finally	
		\begin{equation*}
			\frac{\partial^{\alpha} g}{\partial |x|^{\alpha}} = 	\frac{- 1}{2 \Gamma (1 - \alpha) \cos [(1 - \alpha) \frac{\pi}{2}]} \frac{\partial}{\partial x} \int\limits_{- \infty}^{\infty} \frac{ - sgn (x) \cdot \frac{\partial g}{\partial x} (x^{\prime}) d x^{\prime}}{|x - x^{\prime}|^{\alpha}}.
		\end{equation*}
		
	\end{itemize}
	
	In the case $\alpha = 1/2$ the fractional integrals and derivatives are also called {\em semi-integrals}	and {\em semi-derivatives} \cite{sok02}.

	\subsection{Caputo Fractional Derivative} 
	\label{capu}
	The fractional derivatives in the Caputo sense on the left ($_CD_{a+}^{\alpha}$) and on the right ($_CD_{b-}^{\alpha}$) are defined via the RL fractional integral \cite{gri13}
	$(_CD_{a+}^{\alpha} f) = (J_{a+}^{n -\alpha} f^{(n)}) (x)$ and $(- 1)^n (_CD_{b-}^{\alpha} f) = (J_{b-}^{n -\alpha} f^{(n)}) (x)$. It was introduced independently in 1948 by M. Caputo and by A.N. Gerasimov \cite{ger}; later by Dzherbashyan \& Nersesian \cite{dzh}.
	
	The major difference of the Caputo fractional derivative is that the derivative act first on the function and after the integral is evaluated while in the RL approach the derivative act on the integral.
	
	The Caputo fractional derivative is defined as \cite{pod98}
	\begin{equation}
		\label{cd}
		D_{\star}^{\alpha} f (t) = \left \{
		\begin{array}{rl}
			\displaystyle \frac{1}{\Gamma (1 - \alpha)} \int\limits_0^t  (t - x)^{- \alpha} \frac{df(x)}{dt} dt,  & \quad \alpha \in (0,1)\\
			\displaystyle \frac{df (x)}{dt},
			&  \quad \alpha = 1
		\end{array} \right.
	\end{equation}

	The definition of the Caputo derivative (\ref{cd}) is more restrictive than of the RL one
	(\ref{lrld}, \ref{rrld}) since it requires the absolute integrability of the derivative $df (x)/ dt$ \cite{gor97a}.
	
	The Caputo fractional derivative can be considered as the  regularization in the time origin for the RL derivative \cite{gor03}
	\begin{equation*}
		D_{\star}^{\alpha} f (t) = D^{\alpha} f (t) - f (0^+) \frac{t^{- \alpha}}{\Gamma (1 - \alpha)}	
	\end{equation*}
	and satisfies the property of being zero when applied to a constant. 
	
	Yuan \& Agrawal and Singh \& Chatterjee suggested the  alternative definitions of the Caputo fractional derivative \cite{die10}. The first approach is based on the introduction of the auxiliary bivariate function
	$\phi: (0, \infty) \times [a, b] \rightarrow \mathcal{R}$ as
	\begin{equation*}
		\phi (w, x) = (- 1)^{\lfloor n \rfloor} \frac{2 \sin \pi n}{\pi} w^{2 n - 2 \lceil n \rceil + 1} \int_a^x f^{(\lceil n \rceil)} (\tau) e^{- (x - \tau) w^2} d \tau
	\end{equation*}
	where $\lceil \quad  \rceil$ denote the ceiling function, and, finally
	\begin{equation*}
		D_{\star a}^n f (x) = \int\limits_0^{\infty} \phi (w, x) d w.
	\end{equation*}
	
	The second approach is based on expressing the fractional derivative of the the given function in the form of the integral over $(0, \infty)$ with the integrand that can be obtained as the solution of the first-order initial value problem
	\begin{equation*}
		\frac{\partial \phi^{\star} (w, x)}{\partial x} = - w^{\frac{1}{n - \lceil n \rceil - 1}} \phi^{\star} (w, x) + \frac{(- 1)^{\lfloor n \rfloor} \sin \pi n}{\pi (n - \lceil n \rceil - 1)} f^{(\lceil n \rceil)} (x)
	\end{equation*}
	with the initial condition $\phi^{\star} (w, a) = 0$. 
	Thus
	\begin{equation*}
		\phi^{\star} (w, x) = \frac{(- 1)^{\lfloor n \rfloor} \sin \pi n}{\pi (n - \lceil n \rceil - 1)} \int\limits_0^x  f^{(\lceil n \rceil)} (\tau) \exp (- (x - \tau)w^{\frac{1}{n - \lceil n \rceil - 1}} ) d \tau	
	\end{equation*}
	\begin{equation*}
		D_{\star a}^n f (x) =
		\int\limits_0^{\infty} \phi^{\star} (w, x) d w.
	\end{equation*} 
	
	\subsection{Matrix Approach}
	Operations of the fractional calculus can be expressed by matrix \cite{pod00,pod09}. E.g., the left-sided RL or Caputo derivative can be approximated in the nodes in the equidistant discretization net with the  help of the upper triangular strip matrix $B_n^{(\alpha)}$ as \cite{pod00}
	\begin{equation*}
		\left[v_n^{(\alpha)} \; v_{n - 1}^{(\alpha)} \; \dots \; v_1^{(\alpha)} \; v_0^{(\alpha)} \right]^T = B_n^{(\alpha)} \left[v_n \; v_{n - 1} \; \dots \; v_1 \; v_0 \right]^T 
	\end{equation*}
	where
	\begin{equation*}
		B_n^{(\alpha)}	= \frac{1}{\tau^{\alpha}} \begin{bmatrix}
			\omega_0^{(\alpha)}& \omega_0^{(\alpha)}& \dots& \dots& \omega_{n - 1}^{(\alpha)}& \omega_{n}^{(\alpha)}&\\
			0& \omega_0^{(\alpha)}& \omega_0^{(\alpha)}& \dots& \dots& \omega_{n - 1}^{(\alpha)}& \\
			\dots\\
			\dots\\
			0& 0& 0& \dots& \dots&  \omega_{0}^{(\alpha)}&\\
			
		\end{bmatrix}
	\end{equation*}
	Similarily, the right-hand RL or Caputo fractional derivative can be approximated  with the  help of the corresponding lower triangular strip matrix.
	
	\subsection{Caputo \& Fabrizio Fractional Derivatives} 
	Caputo \& Fabrizio \cite{cf1,cf2} proposed the  fractional derivatives without the singular kernel \cite{los15} by replacing the kernel $(t - \tau)^{- \alpha}$ with the  function $\exp (- \alpha/(1 - \alpha))$ that does not have singularity for $t = \tau$ in the definition of the Caputo derivative and replacing the factor $1/\Gamma (1 - \alpha)$ with $M (\alpha)/(1 - \alpha)$.
	
	E.g., the fractional time derivative for $\alpha \in [0, 1]$ and function $f \in L^1 (-\infty,b)$ is
	\begin{equation*}
		\mathcal{D}_t^{\alpha} f (t) = \frac{\alpha M (\alpha)}{1 - \alpha} \int\limits_{- \infty}^t (f (t) - f (\tau)) \exp \left[- \frac{\alpha (t - \tau)}{1 - \alpha} \right] d \tau
	\end{equation*}
	where $M (\alpha)$ is a normalization function such as $M (0) = M (1) =1 $.
	
	\subsection{GC \& GRL derivatives}
	Zhao \& Luo \cite{zha19a} suggested to divide the fractional derivative with different --- singular and non-singular --- kernels (e.g., RL, Caputo, Caputo-Fabrizio, Atangana-Baleanu \cite{at16}\footnote{
		The equation with the Atangana-Baleanu operator is related to the derivatives of distributed order \cite{tat17}.
	} with the kernel
	\begin{equation*}
		k (x, \alpha) = E_{\alpha} \left(- \frac{\alpha}{1 - \alpha} x  \right),
	\end{equation*}
	Atangana-Gomez \cite{at17} with the kernel
	\begin{equation*}
		k (x, \alpha) = \exp  \left(- \frac{\alpha}{1 - \alpha} x^2  \right)
	\end{equation*}
	derivative with the stretched exponential kernel \cite{sun17} (that is useful in the study of the water diffusion in the human brain using the magnetic resonance imaging \cite{ben06}) 
	\begin{equation*}
		k (x, \alpha) = \exp  \left(- \frac{\alpha}{1 - \alpha} x^{\beta}  \right), \quad \beta > 0, \quad \beta \ne 1)
	\end{equation*} 
	into two classes --- GC (general, Caputo sense) and GRL (general, RL) derivatives that obeys the the principles formulated by V. Volterra in his "general laws of heredity" \cite{volt}:  
	{the linearity principle,}
	{the invariance principle,}
	{the fading memory principle,}
	{the compatibility principle.}
		The compatibility principle requires the validity of two limits: $D_{\alpha} f (x) \rightarrow f (x)$ when $\alpha \rightarrow 0$ and $D_{\alpha} f (x) \rightarrow f^{\prime} (x)$ when $\alpha \rightarrow 1$.
	
	The principle of {\em nonlocality} was suggested by Tarasov \cite{tar18}.
	
	\subsubsection{GC derivatives}
	Zhao \& Luo \cite{zha19a} defined the GC derivative by
	\begin{equation*}
		D_{a, \alpha}^{GC} f (x) = N (\alpha) \int\limits_a^x k (x - t, \alpha)\frac{d f (t)}{d t} d t
	\end{equation*}
	The fading memory principle requires that the remote time and position has less effect: $k (x - t, \alpha)$ decreases when $x$ increases and $k (x - t, \alpha) \rightarrow 0$ when $x \rightarrow \infty$. 
	
	The compatibility principle requires that $N (\alpha) / k (x, \alpha) \rightarrow 1 $ when $\alpha \rightarrow 0$ and $N (\alpha) / k (x, \alpha) \rightarrow \delta(x) $ when $\alpha \rightarrow 1$.
	
	\subsubsection{GRL derivatives}		
	\begin{equation*}
		D_{a, \alpha}^{RL} f (x) = \frac{d}{d x} N (\alpha) \int\limits_a^x k (x - t, \alpha) f (t) d t.
	\end{equation*}	
	
	The restrictions on $k (x - t, \alpha)$ and $N (\alpha)$ are the same.	
	
	\subsection{Marchaud-Hadamard Fractional Derivatives} 
	Marchaud's approach is based on the analytic coninuation of the fractional integrals to the negative orders using the  Hadamard's finite parts of the divergent integrals (Hadamard's idea is to ignore the unbounded contribution to the integral and to assign the value of the remaining --- finite --- expression \cite{die10}).
	
	The Marchaud fractional derivative with the lower limit $a$ is
	\begin{equation*}
		(M_{a+}^{\alpha} f) (x) = \frac{f (x)}{\Gamma (1 - \alpha) (x - a)^{\alpha}} + \frac{\alpha}{\Gamma (1 - \alpha)} \int\limits_a^x \frac{f (x) - f(y)}{(x - y)^{\alpha + 1}} d y
	\end{equation*}
	and with the upper limit $b$ is
	\begin{equation*}
		(M_{b-}^{\alpha} f) (x) = \frac{f (x)}{\Gamma (1 - \alpha) (b - x)^{\alpha}} + \frac{\alpha}{\Gamma (1 - \alpha)} \int\limits_x^b \frac{f (x) - f(y)}{(x - y)^{\alpha + 1}} d y.
	\end{equation*}
	
	Marchard's method is to extend the RL integral  to $\alpha < 0$
	\begin{equation}
		\label{283}
		(J_+^{- \alpha} f) (x) = \frac{1}{\Gamma (- \alpha)} \int\limits_0^{\infty} y^{- \alpha - 1} f (x - y) d y
	\end{equation}
	and to substract the divergent part of the integral in (\ref{283})
	\begin{equation*}
		\int\limits_{\epsilon}^{\infty} y^{- \alpha - 1} f (x - y) d y	= \frac {f (x)}{\alpha \epsilon^{\epsilon}}
	\end{equation*}
	to get finally
	\begin{equation}
		\label{285}
		(M_{+}^{\alpha} f) (x) = \lim_{\epsilon \rightarrow 0 +} \frac{1)}{\Gamma ( - \alpha)}  \int\limits_{\epsilon}^{\infty} \frac{f (x) - f(y)}{y^{\alpha + 1}} d y.	
	\end{equation}
	There are two approaches to extend the definition (\ref{285}) to the case $\alpha > 1$ \cite{hil08}:
	\begin{enumerate}
		\item {To  apply (\ref{285}) to the $n$th derivative $d^nf/d x^n$ for $n < \alpha < n + 1$.}
		\item {To consider $f (x - y) - f (x)$ as the first-order difference and to generalize to $n$th order difference (difference quotient)}
		\begin{equation}
			\label{dq} 
			(\Delta_h^n f) (x) = (\mathcal{I} - Th)^n f (x) = \sum_{k = 0}^n (- 1)^k \binom{n}{k} f (x - k h)
		\end{equation}
		where $\mathcal{I}$ is the identity operator and $T_h = f (x - h)$ is the translation operator.
		
		Thus Marchard fractional derivative for $0 < \alpha < n$ is written as
		\begin{equation*}
			(M_{+}^{\alpha} f) (x) = \lim_{\epsilon \rightarrow 0 +} \frac{1}{C_{\alpha,n}}  \int\limits_{\epsilon}^{\infty} \frac{\Delta_y^n f  (x)}{y^{\alpha + 1}} d y	
		\end{equation*}
		where
		\begin{equation*}
			C_{\alpha,n} = \int\limits_0^{\infty} \frac{(1 - e^{- y})^n}{y^{\alpha + 1}}.	
		\end{equation*}
	\end{enumerate}
	
	\subsection{Gr\"unwald - Letnikov  Derivative} 
	The approach suggested independently by Gr\"unwald in 1867 and Letnikov \cite{let} in 1868 is based on the use the limits of the difference quotients (\ref{dq}) similar to the classical definition of the derivatives for $n \in \mathcal{N}, \quad f \in C^n [a, b], \quad a < x \le b$
	\begin{equation*}
		\tilde{D}^n f (x) = \lim_{h \rightarrow 0} \frac{\Delta_h^n f) (x)}{h^n}
	\end{equation*}
	and extension to the case of the arbitrary $n$.
	
	Since $\binom{n}{k} = 0$ if $n \in \mathcal{N}$  and $n < k$   the expression    (\ref{dq}) is equivalent to 
	\begin{equation}
		\label{dq1} 
		(\Delta_h^n f) (x) = \sum_{k = 0}^{\infty} (- 1)^k \binom{n}{k} f (x - k h).
	\end{equation}
	The series  (\ref{dq1}) is uniformly convergent for any bounded function if  $n > 0$ \cite{gol10}.

	The use of (\ref{dq1}) introduce two problems \cite{die10}: {the function $f$ needs to be defined on $(\infty,b]$; the function $f$ should be such that the series converges.
	
	These problems are solved by the introduction a new function $f^{\star}$
	\begin{equation*}
		f^{\star} = \begin{cases}
			f (x) & x \in [a,b]\\
			0 & x \in (- \infty, a)
		\end{cases}
	\end{equation*}
	that is used instead of the original $f$.
	
	It is also assumed that in the tending to zero $h$ takes only  the Gr\"unwald-Letnikov fractional derivative of order $n$  defined as 
	\begin{equation}
		\label{ll}
		\tilde{D}_{a}^n = \lim_{N \rightarrow \infty} \frac{\Delta_{h_N}^n f (x)}{h_N^n} = \lim_{N \rightarrow \infty} \sum_{k = 0}^N (- 1)^k \binom{n}{k} f (x - k h_N).
	\end{equation}
	
	The Gr\"unwald-Letnikov derivative is called {\em pointwise} or {\em strong} depending on whether the limit is taken pointwise or in  the norm of a suitable Banach space \cite{hil08}.
	
	The binomial coefficient can be generalized to the non-integer arguments
	\begin{equation*}
		(- 1)^j \binom{q}{j} = (- 1)^j \frac{\Gamma (q + 1)}{\Gamma (j + 1) \Gamma (q - j + 1)} = \frac{\Gamma (j - q)}{\Gamma (- q) \Gamma (j + 1)}.
	\end{equation*}
	
	The (left-sided) Gr\"unwald-Letnikov derivative could be written as ($nh = x - a$)
	\begin{equation*}
		\tilde{D}_{a}^{\alpha} f (x) = \lim_{h \rightarrow 0} \frac{1}{h^{\alpha}} \sum_{k =  0}^{\lfloor n \rfloor} (- 1)^k \frac{\Gamma (\alpha + 1) f (x - k h)}{\Gamma (k + 1) \Gamma (\alpha - k + 1)}	
	\end{equation*}
	and right-sided ($nh = b - x$) as
	\begin{equation*}
		\tilde{D}_{b}^{\alpha} f (x) =	 \lim_{h \rightarrow 0} \frac{1}{h^{\alpha}} \sum_{k =  0}^{\lfloor n \rfloor} (- 1)^k \frac{\Gamma (\alpha + 1) f (x + k h)}{\Gamma (k + 1) \Gamma (\alpha - k + 1)}.
	\end{equation*}

	The Gr\"unwald-Letnikov integral of the order $n$ of the function $f$ is written as
	\begin{equation*}
		\tilde{J_a^n} f (x) = \frac{1}{\Gamma (n)} \lim_{N \rightarrow \infty} h_N^n \sum_{k = 0}^N \frac{\Gamma (n + k)}{\Gamma (k + 1)} f (x - k h_N).
	\end{equation*}
	
	\subsection{Riesz  Fractional Operators}
	The fractional integral of the order $\alpha$ in the Riesz sense (also known as the Riesz potential) is defined by the Fourier convolution product
	\begin{equation*} 
		(\mathcal{I}^{\alpha} f) (\vec{x}) = \int\limits_{\mathcal{R}^n} \vec{K}_{\alpha} (\vec{x} - \vec{\xi}) f (\vec{\xi}) d \vec{\xi},
	\end{equation*}
	where $Re (\alpha) > 0$.
	The Riesz kernel
	\begin{equation*}
		\vec{K}_{\alpha} = \frac{1}{\gamma_n (\alpha)}
		\begin{cases}
			\Vert{\vec{x}}\Vert^{\alpha - n}, & \alpha - n \ne 0, 2, \dots,\\	
			\Vert{\vec{x}}\Vert^{\alpha - n} \ln \left(\frac{1}{\Vert{\vec{x}\Vert}} \right), & \alpha - n = 0, 2, \dots
		\end{cases}	
	\end{equation*}
	where $\gamma_n (\alpha)$ is defined by
	\begin{equation*}
		\frac{\gamma_n (\alpha)}{2^{\alpha} \pi^{\frac{n}{2}} \Gamma(\alpha/2)} = 
		\begin{cases}
			\left[\Gamma \left(\frac{n - \alpha}{2}  \right)  \right]^{- 1},	& \alpha - n \ne 0, 2, \dots,\\
			(- 1)^{\frac{n - \alpha}{2}} 2^{- 1}\Gamma \left(\frac{\alpha - n}{2}  \right), & \alpha - n = 0, 2, \dots.
		\end{cases}		
	\end{equation*}
	
	The Riesz fractional integral is \cite{gri13}
	\begin{equation*}
		(\mathcal{I}^{\alpha} f) (\vec{x}) = \frac{\Gamma \left(\frac{1 - \alpha}{2}  \right) }{2^{\alpha} \pi^{\frac{n}{2}}\Gamma(\alpha/2)}	\int\limits_{- \infty}^{\infty} f (\xi) |x - \xi|^{\alpha - 1} d \xi. 
	\end{equation*}
	
	The Riesz fractional derivative is \cite{oli14}
	\begin{equation*}
		D^ {\alpha} [f (x)] = - \frac{1}{2 \cos (\alpha \pi/2)} \frac{1}{\Gamma (\alpha)}
	\end{equation*}	
	\begin{equation*}	
		\frac{d^n}{d x^n} \left[\int\limits_{- \infty}^x (x - \xi)^{n - \alpha_n - 1} f (\xi) d \xi + \int\limits_x^{\infty} (x - \xi)^{n - \alpha_n - 1} f (\xi) d \xi \right].
	\end{equation*}
	
	The Riesz derivative is the generalization of the Laplace operator \cite{zas08} 
	\begin{equation*}
		(- \Delta)^{\frac{\alpha}{2}} =   - \frac{1}{2 \cos (\alpha \pi/2)} \left[ \frac{d^{\alpha}}{d x^{\alpha}} + \frac{d^{\alpha}}{d (- x)^{\alpha}} \right], \quad \alpha \ne 1.
	\end{equation*}

	The Riesz derivative  could be expressed in terms of the Marchaud derivative
	\begin{equation*}
		D^ {\alpha} [f (x)] = - \frac{1}{2 \cos (\alpha \pi/2)} [(M_+^{\alpha} f) (x) + (M_-^{\alpha} f)].	
	\end{equation*}
	
	The related  Riesz-Feller derivative \cite{gor02a} has an additional parameter - "skewness" $\theta$
	\begin{multline*}
		D_{\theta}^{\alpha} f (x) = \frac{\Gamma (1 + \alpha)}{\pi} \times\\ \left[\sin \left[(\alpha + \theta) \frac{\pi}{2} \right]  \int\limits_0^{\infty} \frac{f (x + \xi) f (x)}{\xi^{1 + \alpha}} d \xi 
		+ \sin \left[(\alpha - \theta) \frac{\pi}{2} \right]  \int\limits_0^{\infty} \frac{f (x - \xi) f (x)}{\xi^{1 + \alpha}} d \xi  \right].	
	\end{multline*}
	
	The allowed region of the parameters $\alpha$ and $\theta$ turn out to be a diamond in the plane $\{\alpha, \theta \}$ with the vertices in the points (0,0), (1,1), (2,0), (1, -1) called the "Feller-Takayasu diamond" \cite{gor03,met04}.  
	
	\subsection{Weyl  Fractional Derivative}
	The Weyl  derivative is based on the generalization of the differentiation in the Fourier space \cite{gol10}  --- the integer derivative of the $n$th order $(ik)^n$ of the absolutely integrable function on $[- \pi,\pi]$ presented as the Fourier series is extended to the noninteger  $n$.

	The Weyl  fractional derivative is defined as \cite{mai04}
	\begin{equation*}
		D_{\pm}^{\alpha} = \begin{cases}
			\displaystyle{
				\pm\frac{d}{d x} [I_{\pm}^{1 - \alpha} f (x)]} &0 < \alpha < 1,\\
			
			\displaystyle{
				\frac{d^2}{d x^2} [I_{\pm}^{2 - \alpha} f (x)]} &1 < \alpha < 2,
		\end{cases}		
	\end{equation*}
	where the Weyl fractional integrals are ($\mu > 0$)
	\begin{equation*}
		I_+^{\mu} = \frac{1}{\Gamma (\mu)} \int\limits_{- \infty}^x (x - \chi)^{\mu - 1} f (\chi)d \chi.
	\end{equation*}
	
	\subsection{Erd\'elye-Kober  Fractional Operators}
	The Erd\'elye-Kober integral for a well-behaved function $\phi (t)$ is defined as \cite{luch,pag12}
	\begin{equation*}
		I_{\eta}^{\gamma, \mu} \phi (t)	= \frac{\eta}{\Gamma (\mu)}  t^{- \eta (\mu + \gamma)} \int\limits_0^t \tau^{\eta (\gamma + 1) - 1} (t^{\eta} - \tau^{\eta})^{\mu - 1} \phi (\tau) d \tau,
	\end{equation*}
	where $\mu > 0, \quad \eta > 0, \quad \gamma \in \mathcal{R}$.
	
	In the special case $\gamma = 0, \eta = 1$ the Erd\'elye-Kober  fractional integral is related to the RL fractional integral of the order $\mu$ as
	\begin{equation*}
		I_1^{0, \mu} \phi (t) = \frac{t^{- \mu}}{\Gamma (\mu)} \int\limits_0^t (t - \tau)^{\mu - 1} \phi (\tau) d \tau = t^{- \mu} J^{\mu} \phi (t).
	\end{equation*}
	
	The Erd\'elye-Kober  fractional derivative for $n - 1 < \mu <n, \quad n \in \mathcal{N}$ is defined as
	\begin{equation*}
		D_{\eta}^{\gamma, \mu} \phi (t)	= \prod_{j = 1}^n \left( \gamma + j + \frac{1}{\eta} t \frac{d}{d t} \right) (I_{\eta}^{\gamma + \mu, n - \mu} \phi (t)).	
	\end{equation*}
	
	The Erd\'elye-Kober  fractional derivative reduces to the identity operator when $\mu = 0$
	\begin{equation*}
		D_{\eta}^{\gamma, 0} \phi (t)	= \phi (t)	
	\end{equation*}
	and for $\eta = 1$ and $\gamma = - \mu$ is related to the RL fractional  derivative as 
	\begin{equation*}
		D_{\eta}^{\gamma, \mu} \phi (t)	= t^{\mu} D_{RL}^{\mu} \phi (t).
	\end{equation*}

	\subsection{Interpretation of  Fractional Integral and  Derivatives} 
	
	The integer-order  and integrals have a clear physiscal and geometrical interpretation that simplify their use in practice.
	The numerous different interpretations of the fractional derivatives and integrals have been proposed \cite{hil19}: the probabilistic  \cite{t13, t14a, t15}, geometric \cite{t16, t17, t20, t21}, physical interpretations \cite{t20, t21, t22, t23, t24}.
	
	However, as noted Podlubny \cite{pod02}, "since the appearance of the idea of differentiation and  of arbitrary (not necessary integer) order there was not any acceptable  geometric and physical interpretation of these operations for more than 300
	years". 
	
	Teneiro Machado \cite{t15} wrote the G\"unwald-Letnikov derivative of $x (t)$ as
	\begin{equation*}
		D^{\alpha} [x (t)] = \lim_{h \rightarrow 0}  \left[\frac{1}{h^{\alpha}} \sum_{k = 0}^{\infty} \gamma (\alpha, k) x (t - k h) \right], \quad
		\gamma (\alpha, k) = (- 1)^k \frac{\Gamma (\alpha + 1)}{k! \Gamma (\alpha -_{} k + 1)}
	\end{equation*}
	where $h$ is the time increment.
		The author noted that
	\begin{equation*}
		\gamma (\alpha, 0) = 1, \quad - \sum_{k = 0}^{\infty} \gamma (\alpha, k) =1
	\end{equation*}
	thus the "present" (P) is constituted by $x (t)$  the probability 1 while the totality of the "past/future"
	(PF) is constituted by the samples $x (t), x (t - h), x (t - 2h), \dots$; each sample is weighted with a probability $- \gamma (\alpha, k)$. 
	
	Nigmatullin \cite{t22,nig05} interpreted the fractional integral in terms of the fractal Cantor set. The author considered the evolution of the state of the physical system
	\begin{equation*}
		J (t) = \int\limits_0^t K (t, \tau) f (t) d \tau
	\end{equation*}
	where the memory function $K (t, \tau) f (t)$ accounts for the loss of some states of th system; the fractional index of integration equals the fractal dimension of the Cantor set.

	Podlubly \cite{pod02} and Podlubny et al. \cite{pod07} suggested the geometrical interpreation of the left-sided (equation (\ref{lrl})) and right-sided (equation (\ref{rrl})) RL integrals and of the RL (equations (\ref{lrld}) - (\ref{rrld})) and the Caputo (equation (\ref{cd})) derivatives, as well as of the Riesz potential that is the sum of the left-sided and right-sided RL fractional integrals
	\begin{equation*}
		R_b^{\alpha} f (x) = \frac{1}{\Gamma (\alpha)} \left[ \int\limits_a^x (x - t)^{\alpha - 1} f (t) d t + \int\limits_x^b (t - x)^{\alpha - 1} f (t) d t \right]
	\end{equation*}
	and of the Feller potential
	\begin{equation*}
		\Phi^{\alpha} f (x) = c J_{a^+}^{\alpha} f (x) + 	d J_{b^-}^{\alpha} f (x).
	\end{equation*}
	
	The geometric interpretation by Podlubny is based on additing the third dimension
	\begin{equation*}
		g_x (t) = \frac{1}{\Gamma (\alpha +1)} [x^{\alpha} + (x - t)^{\alpha}]
	\end{equation*}
	to the pair (t, f (x)) and considering the three-dimensional line $(t, g_x (t), f (t))$ as the top edge of the "fence" that gives shadow on the wall in the (g,f) plane. 
	
	Tarasov \cite{tar17}  proposed the ”informatic” (”computer science”) interpretation of the RL and the Caputo derivatives of the non-integer orders using the reconstructions from the infinite sequence of the derivatives of the integer orders. Such reconstructions atre based on the  Kotel’nikov theorem (also known as the sampling theorem)  proved by Vladimir Kotel’nikov  in 1933 and also  by Claude Shannon 1949:
	under the certain restrictive conditions, function $f(t)$ can be restored from its samples $f[n] = f(nT)$ according to the Whittaker-Shannon interpolation formula. The author stressed that infinity of the sequences of the integer derivatives plays a fundamental role in representation of the fractional derivatives that describe nonlocality and memory.
	
	G\'omez-Aguilar et al. \cite{gom14} analysed the Caputo differentiation using the RC circuit for which the fractional version of the Ohm's law and Kirchhoff's law are written as
	\begin{equation*}
		v (t) = \frac{1}{\sigma^{1 - \gamma}} \frac{d^{\gamma} q}{d t^{\gamma}}, \qquad R \frac{d q}{d t} + \frac{1}{C} q (t) = v (t)
	\end{equation*}
	where $q$ is the electric charge, $v$ is the voltage, $R$ is the resistance of the conductor, $C$ is the capacitance. The parameter $\sigma$ is introduced in order to be consistent with the dimensionality; it characterizes the fractional structures (the components that show the intermediate behaviour between conservative (capacitor) and dissipative (resistor) \cite{gom14}. The authors derived the fractional differential equation for the RC circuit
	\begin{equation*}
		\frac{d^{\gamma}}{d t^{\gamma}} + \frac{1}{\tau_{\gamma}} q (t) = \frac{C}{\tau_{\gamma}} v (t), \qquad
		\tau_{\gamma} = \frac{R C}{\sigma^{1 - \gamma}}	
	\end{equation*}
	where $\tau_{\gamma}$ is the time constant. 
	G\'omez-Aguilar et al. claimed that the  differentiation is related to the  memory effects that reflect the intrinsic dissipation in the system.
	
	Sierociuk et al. \cite{sie15} used the RC network to model the fractional order diffusion based on the analogy between the heat and electrical conduction. The authors 
	showed that the equations for the capacitor and for the resistor in the transmission line could be used to get the diffusion equation; the loosing of heat was represented by the additional resistors connected parallel to capacitors.

	Carpinteri et al. \cite{car09} considered the mechanical interpretation of the Marchaud fractional derivative using the body springs connecting the non-adjacent points of the body with the stiffness decaying with the distance between the material points.

	\subsection{Local Fractional Derivatives} 
	The fractional derivatives are nonlocal. Several researches introduced the local variants \cite{zha06} that are useful for study of the pointwise behaviour of the fractal and multifractal functions that describe, e.g., the stress and deformation patterns in materials exhibiting the fractal-like microstructure \cite{car09} or the velocity field of turbulent fluid \cite{kol98}.  
	
	Kolwankar \& Gangal \cite{kol97,kol98,kol98a,kol13} defined
	the derivative via the RL derivative as
	\begin{equation*}
		\mathfrak{D}^q f (y) = \lim_{x \rightarrow y} \frac{D^q (f (x) - f (y))}{(x - y)^q}
	\end{equation*}
	if the limit exists and  finite.
	
	The local fractional Taylor formula is written as \cite{yang12c}  
	\begin{equation*}
		f (x) = \sum_{i = 0}^n \frac{f^{(n) (y)}}{\Gamma (1 + n)} (x - y)^n  \frac{\mathfrak{D}^{\alpha}}{\Gamma (n! + \alpha)} (x - y)^{\alpha} + R_{\alpha} (x, y).
	\end{equation*}
	
	Yang et al. \cite{yang16,zha16} used similar definition
	\begin{equation*}
		\mathfrak{D}^{(k)} f (\tau) = \lim_{\tau \rightarrow \tau_0}  \frac{f (\tau) - f (\tau_0)}{\tau^k - \tau_0^k}.	
	\end{equation*}

	Chen et al. \cite{che10} proposed the local derivatives based on the integrals of the difference-quotient (IDQ) or the singular  of difference-quotient (SIDQ). For example,  the right SIDQ local derivative is
	\begin{equation*}
		\mathfrak{D}^{\alpha} f (y) = \frac{1}{\Gamma (1 - \alpha)} \lim_{h \rightarrow 0_+} \int\limits_0^1 (1 - t)^{- \alpha} \frac{f (t h + y) - f (y)}{h^{\alpha}} d t. 
	\end{equation*}
	
	The local fractional derivative is essentially the {\em fractal} derivative \cite{chen06a,he14}. In contrast to the purely analytical approach of the fractional calculus, the fractal calculus follows the physical-geometric approach; to avoid confusion it is suggested to call the latter the {\em scaled} calculus \cite{ont13}.
	
	The fractal ("Hausdorff") derivative on the time fractal is defined as \cite{he18}
	\begin{equation*}
		\frac{\partial f}{\partial t^{\sigma}} = \lim_{t_B \rightarrow t_A} \frac{f (t_B) - f (t_A)}{(t_B)^{\sigma} - (t_A)^{\sigma}}
	\end{equation*}
	where ${\sigma}$ is the fractal dimension of time.
	
	A more general definition is formulated as \cite{h1,h2}
	\begin{equation*}
		\frac{\partial^{\tau} f}{\partial t^{\sigma}} = \lim_{t_B \rightarrow t_A} \frac{f^{\tau} (t_B) - f^{\tau} (t_A)}{(t_B)^{\sigma} - (t_A)^{\sigma}}
	\end{equation*}

	Since the fractal derivative is the local operator, the  numerical solution of the fractal derivative equations can be performed by the standard numerical techniques for the integer-order differential equations \cite{chen10}. 
	The similar properties have the so called "conformable"  fractional derivatives.
	
	\subsubsection{"Conformable"  Fractional Derivative}
	Most fractional derivatives do not have the desirable properties \cite{abd14,kha14,kat14}:
	 {the derivative of a constant is not zero;}
		{they do not obey the product rule $D^{\alpha} (f g) = f D^{\alpha} (g) + g D^{\alpha} f$ ; }
		{they do not obey the quotient rule
		$D^{\alpha} ({f}/{g}) = (g D^{\alpha} (f) - f D^{\alpha} (g)) / g^2$; 		they do not obey 
		the chain rule $D^{\alpha} (f g) = f^{\alpha} (g (t) g^{\alpha} (t)$;
		they do not obey in general $D^{\alpha} D^{\beta} f = D^{\alpha + \beta} f$.
	
	Khalil et al.\cite{kha14} and Katugampola \cite{kat14,kat14a}
	suggested the so called "conformable" limit based \cite{and15} derivatives
	\begin{equation*}
		D^{\alpha} f (t) = \lim_{\epsilon \rightarrow 0} \frac{f (t + \epsilon t^{1 - \alpha}) - f (t)}{\epsilon}, \qquad 0 < \alpha < 1, 	
	\end{equation*}
	and
	\begin{equation*}
		D^{\alpha} f (t) = \lim_{\epsilon \rightarrow 0} \frac{f (t e^{\epsilon t^{-\alpha}}) - f (t)}{\epsilon}, \qquad 0 < \alpha < 1. 		
	\end{equation*}
	
	Since the conformable derivative is the extension of the classical derivative definition, this derivative obeys the product rule, the quotient rule, the linearity property, the zero derivative for the constant and are  valid for some extensions of the classical calculus such as the   Rolle's Theorem or Mean Value Theorem \cite{iyi16}.
	
	\section{Tempered Fractional Calculus}
	Sabzikar et al. \cite{temp} suggested a variant of the fractional calculus where power laws are tempered by the exponential factor. The random walks model with the exponentially tempered power law jumps converges to a tempered stable motion \cite{cha11}. This {\em tempered} fractional diffusion is useful in the geophysical \cite{s44,s70} and financial \cite{s11} problems.
	
	The authors considered two intervals for the parameter $\alpha$:
	\begin{itemize}
		\item {$0 < \alpha < 1$. 			
			The {\em tempered}	fractional derivative $\partial_x^{\alpha, \lambda}$ is defined as the function with the Fourier transform $[(\lambda + i k)^{\alpha} - \lambda^{\alpha}] \hat{f} (k)$ that in real space is written as
			\begin{equation*}
				\partial_x^{\alpha, \lambda} f (x) = \frac{\alpha}{\Gamma (1 - \alpha)}	\int\limits_0^{\infty} (f (x) - f (x - y)) e^{- \lambda y} y^{- \alpha - 1} d y.
			\end{equation*}
			
			The negative {\em tempered}	fractional derivative $\partial_{- x}^{\alpha, \lambda}$ is defined as the function with the Fourier transform $[(\lambda - i k)^{\alpha} - \lambda^{\alpha}] \hat{f} (k)$ that in real space is written as
			\begin{equation*}
				\partial_{- x}^{\alpha, \lambda} f (x) = \frac{\alpha}{\Gamma (1 - \alpha)}	\int\limits_0^{\infty} (f (x) - f (x + y)) e^{- \lambda y} y^{- \alpha - 1} d y.
			\end{equation*}
		}
		\item {$1 < \alpha < 2$. 			
			The {\em tempered}	fractional derivative $\partial_x^{\alpha, \lambda}$ is defined as the function with the Fourier transform $[(\lambda + i k)^{\alpha} - \lambda^{\alpha} - i k \alpha \lambda^{\alpha - 1}] \hat{f} (k)$ that in real space is 
			\begin{equation*}
				\partial_x^{\alpha, \lambda} f (x) = \frac{\alpha (1 - \alpha)}{\Gamma (2 - \alpha)}	\int\limits_0^{\infty} (f (x - y) - f (x) + y f^{\prime} (x)) e^{- \lambda y} y^{- \alpha - 1} d y.
			\end{equation*}
			
			The negative {\em tempered}	fractional derivative $\partial_x^{\alpha, \lambda}$ is defined as the function with the Fourier transform $[(\lambda - i k)^{\alpha} - \lambda^{\alpha} + i k \alpha \lambda^{\alpha - 1}] \hat{f} (k)$ that in real space is 
			\begin{equation*}
				\partial_{- x}^{\alpha, \lambda} f (x) = \frac{\alpha (1 - \alpha)}{\Gamma (2 - \alpha)}	\int\limits_0^{\infty} (f (x + y) - f (x) - y f^{\prime} (x)) e^{- \lambda y} y^{- \alpha - 1} d y.
			\end{equation*}
		}
	\end{itemize}
	
	Sabzikar et al. introduced the positive tempered integral as
	\begin{equation*}
		\mathfrak{I}_+^{\alpha, \lambda} f (x) = \frac{1}{\Gamma (\alpha)} \int\limits_{- \infty}^x f (u) (x - u)^{\alpha - 1}e^{- \lambda (x - u)} d u
	\end{equation*}
	
	and the negative tempered integral as
	\begin{equation*}
		\mathfrak{I}_-^{\alpha, \lambda} f (x) = \frac{1}{\Gamma (\alpha)} \int\limits_{- \infty}^x f (u) (u - x)^{\alpha - 1}e^{- \lambda (u -x)} d u
	\end{equation*}
	called the RL tempered integrals since for $\lambda = 0$ they reduce to the usual RL integrals.
	
	The authors defined the RL tempered fractional derivatives $\mathcal{D}_{\pm}^{\alpha, \lambda}$ as functions with the Fourier transform $(\lambda \pm i k)^{\alpha} \hat{f} (k)$ that can be expressed  
	\begin{equation*}
		\mathcal{D}_{\pm}^{\alpha, \lambda} f (x )= 
		\begin{cases}
			\partial_{\pm x}^{\alpha, \lambda} f (x) + \lambda^{\alpha} f (x) & 0 < \alpha < 1 \\
			\partial_{\pm x}^{\alpha, \lambda} f (x) + \lambda^{\alpha} f (x) \pm \alpha \lambda^{\alpha - 1} f^{\prime} (x)\ &1 < \alpha  < 2.
		\end{cases}
	\end{equation*}
	
	Evidently, integration and differentiation are the inverse operators: $$\mathcal{D}_{\pm}^{\alpha, \lambda} \mathfrak{I}_{\pm}^{\alpha, \lambda} f (x) = f (x), \quad  \mathfrak{I}_{\pm}^{\alpha, \lambda} \mathcal{D}_{\pm}^{\alpha, \lambda} f (x) = f (x).$$
	
	The integration and differentiation operators have the semigroup property
	\begin{equation*}
		\mathfrak{I}_{\pm}^{\alpha, \lambda} \mathfrak{I}_{\pm}^{\beta, \lambda} f = \mathfrak{I}_{\pm}^{\alpha + \beta, \lambda} f, \quad 
		\mathcal{D}_{\pm}^{\alpha, \lambda} \mathcal{D}_{\pm}^{\beta, \lambda} f = \mathcal{D}_{\pm}^{\alpha + \beta, \lambda} f.
	\end{equation*}

	\section{Fractional Differential Equations}
	\index{fractal}\index{fractional}
	Generally, the fractal media could not be considered as continuous media.
	The use of the non-integer dimensional spaces \cite{old74} is necessary to describe a fractal medium by the continuous models \cite{tar16}.
	The fractional differential equations \cite{duan,nakh06,kilb06,die10}  are non-local (i.e. could incorporate the effects of the memory and the spatial correlations) and  could be formulated in the  distinct but mathematically equivalent forms. Mainardi et al. \cite{main07} compared the fractional extensions of the standard Cauchy problem 
	\begin{equation}
		\label{m21}
		\frac{\partial u (x, t)}{\partial t} = \frac{\partial^2 u (x, t)}{\partial x^2}, \qquad x \in \bf{R}, \qquad t \in R_0^+, \qquad 
		u (x, 0^+) = u_0 (x).
	\end{equation}
		
	The fundamental solution (or Green function) of  (\ref{m21}), i.e. the solution subjected to the initial condition $u_0 (x) = \delta (x)$, is the Gaussian probability density function
	\begin{equation*}
		u (x, t) = \frac{1}{2 \sqrt{\pi}} t^{- 1/2}e^{- x^2 / (4 t)}.
	\end{equation*}
	
	The Green function has the scaling property
	$u (x, t) = t^{{1}/{2}} U ({x}/t^{1/2})$,	
	$U (x)$ is the reduced Green function.
	
	The Cauchy problem (\ref{m21}) is equivalent to the integro-differential equation
		\begin{equation*}
			u (x, t) = u_0 (x) + \int\limits_0^t \left[\frac{\partial^2 u (x,\tau)}{\partial x^2} \right] d \tau	
		\end{equation*} 
		where the initial condition is incorporated.
		
		The fractional diffusion equation could be written with the use of the RL derivative  $D^{1 - \beta}$ ($\beta$ is the real number $0 < \beta < 1$)
		\begin{equation}
			\label{m27}
			\frac{\partial u (x, t)}{\partial t} = D^{1 - \beta} \frac{\partial^2 u (x, t)}{\partial x^2}
		\end{equation}
		or the Caputo derivative $D_{\star}^{\beta}$
		\begin{equation}
			\label{m28}
			D_{\star}^{\beta} u (x, t) = \frac{\partial^2 u (x, t)}{\partial x^2}.	
		\end{equation}
		
		The equations (\ref{m27}) and (\ref{m28}) are equivalent to the equation based on the use of the RL fractional integral of the order $\beta$
		\begin{equation}
			\label{m29}
			u (x, t) = u_0 (x) + J^{\beta} \left[\frac{\partial^2 u (x,\tau)}{\partial x^2} \right].
		\end{equation}
		
		Note that the equation (\ref{m27}) could be obtained by differentiating (\ref{m29}), the equation (\ref{m29}) can be derived by the fractional integration of (\ref{m28}).
		
		The equation (\ref{m27}) was studied by Metzler et al. \cite{met94} and by Saichev \& Zaslavsky \cite{sai97}, the equation (\ref{m29}) by Gorenflo et al. \cite{gor95,gor98} and by Mainardi \cite{m96,m97},
		the integrodifferential equation (\ref{m29}) by Schneider \& Wyss \cite{sch89} using the Mellin transform.
		
		Mainardi et al. \cite{main07} search for  the  fundamental solution 
		of the equation (\ref{m28}) by applying the sequence of the 
		Fourier 
		\begin{equation*}
			\mathcal{F} \{v (x); k\}	= \hat{v} (k) = \int\limits_{- \infty}^{\infty} e^{i k x} v (x) d x, \qquad k \in \bf{R}
		\end{equation*}
		and the Laplace  
		\begin{equation*} 
			\mathcal{L}	\{w (t); s \} =  \tilde{w} (s) = \int\limits_0^{\infty} e^{- s t} w (t) d t, \qquad s \in \bf{C}
		\end{equation*}
		transforms.
				Thus the Green function in the Fourier-Laplace domain is determined by
		\begin{equation}
			\label{m212}
			\hat{\tilde{u}} (k, s) = \frac{s^{\beta - 1}}{s^{\beta} + k^2}, \qquad 0 < \beta \le 1, \qquad \mathcal {R} (s) > 0, \qquad k \in \bf{R}.	
		\end{equation}
		
		There are two strategies to determine the Green function in the space-time domain $u (x,t)$ related to the order in performing inversions in the expression (\ref{m212}) \cite{main07}:
		1) {Invert the Fourier transform to get $\tilde{u} (x, s)$ and then invert the Laplace transform \cite{m96,m97}} or 2) invert the Laplace transform to get $\hat{u} (k,t)$ and then invert the Fourier transform \cite{gor00,main01}.
		
		Nieto \cite{nie10} studied the linear fractional differential equation with the spatial RL derivative for initial or periodic boundary conditions and derived the maximum principle using the properties of the Mittag-Leffler functions. 
				Compte \cite{com96} and West et al. \cite{wes97} studied  
		the equation for the hyperdiffusion (L\'evy-flight diffusion)
		\begin{equation*}
			\frac{\partial P}{\partial t} = D (- \Delta)^{\frac{\gamma}{2}}	
		\end{equation*}
		where the fractional $n$-dimensional Laplace operator $(- \Delta)^{\frac{\gamma}{2}}$ is defined by its Fourier transform with respect to the spatial variable \cite{duan}
		\begin{equation*}
			\mathcal{F}[(- \Delta)^{\frac{\gamma}{2}} g (x)] = |\omega|^{\gamma} \mathcal{F} [g (x)].
		\end{equation*}

		Luchko \cite{luc09} derived the maximum principle for the initial-boundary-value problem for the time-fractional diffusion equation with Caputo derivative over the open bounded domain $G \times (0. T), G \subset R^n$.
		
		The  equation could subjected to the  complex transformation \cite{y16,he12,li12b}
		$s = {x^S}/{\Gamma (1 + \alpha)}$ 
		to convert to a partial differential equation\footnote{
			Such transformation is possible in the multidimensional case if the variables $t, x, y, z$ obey the following constraint \cite{he12} ($q, p, k, l$ are constants)
			\begin{equation*}
				\xi = \frac{q t^{\alpha}}{\Gamma (1 + \alpha)} + \frac{p x^{\beta}}{\Gamma (1 + \beta)} + \frac{k y^{\gamma}}{\Gamma (1 + \gamma)} + \frac{l z^{\lambda}}{\Gamma (1 + \lambda)}. 
			\end{equation*} }
		For example, the heat conduction equation ($\alpha$ is the fractal dimension of the fractal medium)
		\begin{equation*}
			\frac{\partial T}{\partial t} = \frac{\partial^{\alpha}}{\partial x^{\alpha}} \left(\lambda \frac{\partial^{\alpha} T}{\partial x^{\alpha}} \right)
		\end{equation*}
		is converted into the equation
		\begin{equation*}
			\frac{\partial T}{\partial t} = \frac{\partial}{\partial s} \left(\lambda \frac{\partial T}{\partial s} \right).
		\end{equation*}
		That could be further transformed by introduction of the Boltzmann variable \cite{liu15a}
		$\chi = {s}/{\sqrt{t}} = {x^{\alpha}}/{\sqrt{t} \Gamma (1 + \alpha)}$ 
		into the ordinary differential equation
		\begin{equation*}
			\frac{d}{d \chi}  \left(\lambda \frac{d T}{d \chi} \right)	+ \frac{\chi}{2} \frac{d T}{d \chi}.
		\end{equation*}
		
		For the general fractional differential equation in the Jumarie's modification of the RL derivatives $$
			f (u, u_t^{\alpha}, u_x^{\beta}, u_y^{\gamma}, u_z^{\lambda},  u_t^{2 \alpha}, u_x^{2 \beta}, u_y^{2 \gamma}, u_z^{2 \lambda}, \dots) = 0$$
			He \& Li \cite{he13} suggested the following transforms
		\begin{equation*}
			s = \frac{q t^{\alpha}}{\Gamma (1 + \alpha)}, \quad X = \frac{p x^{\beta}}{\Gamma (1 + \beta)}, \quad Y = \frac{k y^{\gamma}}{\Gamma (1 + \gamma)}, \quad Z =  \frac{l z^{\lambda}}{\Gamma (1 + \lambda)} 	
		\end{equation*}
		thus converting the fractional derivatives into classical derivatives
		\begin{equation*}
			\frac{\partial^{\alpha} u}{\partial t^{\alpha}} = q \frac{\partial u}{\partial s}, \quad
			\frac{\partial^{\beta} u}{\partial x^{\beta}} = p \frac{\partial u}{\partial X}, \quad
			\frac{\partial^{\gamma} u}{\partial y^{\gamma}} = k \frac{\partial u}{\partial Y}, \quad
			\frac{\partial^{\lambda} u}{\partial z^{\lambda}} = l \frac{\partial u}{\partial Z}.
		\end{equation*}
		
		Babusci et al. \cite{bab} discussed relations between the differential equations and the theories of the pseudo-operators \cite{po1,po2} and the generalized integral transforms.

		\subsection{Distributed order differential equations}
		\label{dofe_}
		The distributed order differential equations (DODE) are a special class of the fractional differential equations \cite{do1,do2,do3a,do3b,do3,do4,do5,gor05,main07a}.
		Chechkin et al. \cite{chech} discussed the natural and the modifies forms of DODEs and noted that the latter in combination with the continuity equation and the retarded linear response equation for the flux exhibiting memory of the processes at the previous times admits a thermodynamic interpretation. DODEs are used to describe 	the accelerating subdiffusion,   decelerating superdiffusion or transformation of the anomalous behaviour at the short times into the normal behaviour at the long times. For example, Metzler \& Klafter  \cite{met04} considered the DODE for the description of the ultraslow diffusion with the logarithmic time dependence $\left< x^2 (t) \right> \propto \log^k t$ including the so called Sinai diffusion ($k = 4$).
		
		The concept of the distributed order differentiation is close to the variable order fractional operators that are useful for the study of the viscoelasticity, the reaction kinetics of proteins, the electrorheological fluids, the damage modelling \cite{gor02,var1,var2}.
		
		There are two approaches to the formulation of the distributed order differential equations: 
		1)  {direct  --- a new variable does not assigned;} 2) {Independent variable approach --- the order is considered as a function of some independent variable.}
		
		Mainardi et al. \cite{main07} studied the fractional diffusion equation of distributed order
		\begin{equation}
			\label{m31}
			\int\limits_0^1 b (\beta) [D^{\beta} u (x,t)] d \beta = \frac{\partial^2 u (x,t)}{\partial x^2}, \quad b(\beta) \ge 0, \quad \int\limits_0^1 b(\beta) d \beta = 1
		\end{equation}
		with $x \in \bf{R}$, $t \ge 0$ and the initial condition
		$u (x,0^+) = \delta (x)$.
		The weight function $b (\beta)$ is called the order-density. 
			The authors used the Fourier and Laplace transforms to get the fundamental solution  similar to a single-order case (\ref{m212})
		\begin{equation*}
			\left[\int\limits_0^1 b (\beta) s^{\beta} d \beta \right] \hat{\tilde{u}} (k, s) - \int\limits_0^1 b (\beta) s^{\beta - 1} d \beta = - k^2 \hat{\tilde{u}} (k, s)
		\end{equation*}
		and
		\begin{equation}
			\label{m32}
			\hat{\tilde{u}} (k, s) = \frac{B (s)/s}{B (s) + k^2},   \qquad k \in \bf{R}, \qquad	
		B (s) = \int\limits_0^1 b (\beta) s^{\beta} d \beta.
		\end{equation}
		
		In the case of small $k$ the equation (\ref{m32}) can be approximated as
		\begin{equation*}
			\hat{\tilde{u}} (k, s) = \frac{1}{s} \left( 1 - \frac{k^2}{B (s)} +  \dots \right)	
		\end{equation*}
		and the second moment is written as
		\begin{equation}\label{m34}
			\tilde{\mu_2} (s) = - \frac{\partial^2}{\partial k^2}	\hat{\tilde{u}} (k=0, s) = \frac{2}{B (s)}.	
		\end{equation}
		
		The special case of DODEs are the double-order fractional equations \cite{chech}
		\begin{equation*}
			b (\beta) = b_1 \delta(\beta - \beta_1) + b_2 \delta(\beta - \beta_2), \quad 0 < \beta_1 < \beta_2 \le 1, \quad \beta_1 > 0,  \beta_2 >  0, \quad \beta_1 + \beta_2 = 1.
		\end{equation*}
		
		Asymptotic behaviour of $\mu_2 (t)$ follows from (\ref{m34}) for cases of {the slow diffusion (the power-law growth, $b (\beta) = b_1 \delta(\beta - \beta_1) + b_2 \delta(\beta - \beta_2) $) where $	\tilde{\mu_2} (s) =	{2}/({b_1 s^{\beta_1 + 1} + b_2 s^{\beta_2 + 1})}$ and 			}
			 {the ultra-slow diffusion (the logarithmic growth, $b (\beta)  =1$) with 
			$\tilde{\mu_2} (s) =	2 {\ln s}/{s (s - 1)}$.
			}

		The distributed order equations allow to describe the more complex media.
		The time-fractional diffusion equation of the distributed order (\ref{m31}) is potentially more flexible to represent the local phenomena while the space-fractional diffusion equation of the distributive order is more suited to represent the variations in space \cite{cap03}.

		\subsection{Special Functions}
		There are special functions related to the differential equation similar to the classical case (such as e.g., the Bessel and the cylindrical functions, the classical orthogonal polynomials, Airy functions etc.) \cite{leb}. 
		The most important functions in  the fractional calculus are the Mittag-Leffler function \cite{hau11},
		the H-functions \cite{hil0,hf,kir2}, the Wright functions \cite{wf,wf1}, the generalized Lommel-Write functions \cite{lwf}.
		The Mittag-Leffler function  is even called the "Queen"-function of the fractional calculus \cite{kir2}.
		\subsubsection{Mittag-Leffler Functions}
		\label{ml}
		The eigenfunction  of the RL derivatives are the solutions of the equation \cite{hil08}
		\begin{equation*}
			D_{0+}^{\alpha} [ f (x)] = \lambda f (x)	
		\end{equation*}
		where $\lambda$ is the eigenvalue.
		The eigenfunctions are $
			f (x) = x^{1 - \alpha} E_{\alpha, \alpha} (\lambda x^{\alpha})$ 	where
		\begin{equation}
			\label{ml2}
			E_{\alpha, \beta} = \sum_{k = 0}^{\infty} \frac{x^K}{\Gamma (\alpha k + \beta)}
		\end{equation}
		is the generalized  Mittag-Leffler function (also called the Wiman's function \cite{shu07}).
		
		The more general eigenvalue equation for derivatives of the orders $\alpha$ and $\beta$ is
		\begin{equation*}
			D_{0+}^{\alpha, \beta} [ f (x)] = \lambda f (x)		
		\end{equation*}
		The solution is \cite{hil08} $	f (x) = x^{(1 - \beta)(1 - \alpha)} E_{\alpha, \alpha + \beta} (\lambda x^{\alpha})$. 
		
		The special case is the equation $D_{0+}^{\alpha, 1} [ f (x)] = \lambda f (x)$ 
		with eigenfunction $f (x) = E_{\alpha} (\lambda x^{\alpha})$. 
		
		The one-parameter Mittag-Lefler function is the particular case of (\ref{ml2}) for $\beta = \alpha$. 
	
		Evidently \cite{die10},
		\begin{equation*}
			E_{0, 1} (x) = \sum_{k = 0}^{\infty} \frac{1}{\Gamma (1)} = \sum_{k = 0}^{\infty} x^k = \frac{1}{1 - x}, \qquad 	
		E_1 (x) = \sum_{k = 0}^{\infty} \frac{x^k}{k!} = \exp (x). 	
		\end{equation*}
		
		There are other special cases such as \cite{hau11} $E_2 (- x^2) = E_{2, 1} (- x^2) = \cos (x)$; $E_2 (x^2) = E_{2, 1} (x^2) = \cosh (x)$; for $x > 0 \quad E_{1/2} (x^{1/2}) = E_{1/2} (x^{1/2}) = (1 + erf (x) ) \exp (x^2)$;
			for $x \in \mathcal{C}$ and $r \in \mathcal{N}$
				\begin{equation*}
					E_{1, r} = \frac{1}{x^{r - 1}} \left(\exp (x) - \sum_{k = 0}^{r - 2} \frac{x^k}{k!}  \right);
				\end{equation*}	
			$E_3 (x) = {1}/{2} [e^{x^{1/3}} + 2 e^{- 1/2 x^{1/3}} \cos ({\sqrt{3}}/{2} x^{1/3})]$;
			$E_4 (x) = {1}/{2} [\cos(x^{1/4}) + \cosh(x^{1/4})]$ 
		where $	erf (x) = {2}/{\sqrt{\pi}} \int_0^x \exp(-t^2) d t$.
				
		The Mittag-Leffler function $E_1$ satisfies the functional relation \cite{die10,hau11}
		$E_1 (x-y) = {E_1 (x)}/{E_1 (y)}$
		and the relation between two Mittag-Leffler functions with different parameters
		$E_{n_1,n_2} (x) = x E_{n_1, n_1 + n_2} (x) + {1}/{\Gamma (n_2)}$. 
		Note that the frequently used relation $
			E_{\alpha} (a(t + s)^{\alpha}) = E_{\alpha} (a t^{\alpha}) E_{\alpha} (a s^{\alpha}), \quad t, s \ge 1$ 
		is valid only if $\alpha = 0$ or $\alpha = 1$ \cite{pen10}.
		
		Asymptotic expansions and integral representations of the Mittag-Leffler functions could be found in the papers \cite{gor97,gor02,hau11}.
		
		Prabhakar \cite{pra71} suggested the extension
		\begin{equation*}
			E_{\alpha, \beta}^{\gamma} (x) = \sum_{n = 0}^{\infty} \frac{(\gamma)_n}{\Gamma (\alpha n + \beta)} \frac{x^n}{n!}, \qquad Re (\alpha) > 0, Re (\beta) > 0
		\end{equation*}
		$(\gamma)_n$ is the Pochhammer symbol \cite{shu07} $(\gamma)_0 =1, (\gamma)_n = \gamma (\gamma + 1) (\gamma + 2) \dots (\gamma + n - 1)$. 		
		
		The extension to the multi-index Mittag-Leffler functions \cite{kir1,kir2} 
		\begin{equation*}
			E_{(\frac{1}{\rho_i}),(\mu_i)} (x) = \sum_{k = 0}^{\infty} \frac{x^k}{\Gamma (\mu_1 + k/\rho_1) \dots \Gamma (\mu_m + k/\rho_m)}
		\end{equation*}
		is performed by replacing of the indices $\alpha = 1/\rho$ and $\beta = \mu$ by two sets of multi-indices $\alpha \rightarrow (1/\rho_1, 1/\rho_2, \dots, 1/\rho_m)$ and $\beta \rightarrow (\mu_1, \mu_2, \dots, \mu_m)$. 
		
		There are a couple of related functions \cite{nakh03}
		\begin{itemize}
			\item {Barret's function
				\begin{equation*}
					U (x, \lambda) = \sum_{k = 1}^{\infty} \frac{\lambda^{k - 1} x^{k \alpha i}}{\Gamma (k \alpha - i +1);}
				\end{equation*}	
			}
			\item {Rabotnov's (fractional exponential) function \cite{rab97,rab1}}
			\begin{equation*}
				\mathcal{E}_{\alpha} (\beta, x) = x^{\alpha} \sum_{n = 0}^{\infty} \frac{\beta^n x^{n (\alpha + 1)}}{\Gamma ((n + 1)(1 + \alpha))}.
			\end{equation*}
		\end{itemize}

		\subsubsection{H Functions}
		The H-function of order $(m, n, p, q) \in \mathcal{N}^4$  is defined via the Mellin-Barnes type contour integral \cite{duan,hf}
		\begin{equation*}
			H_{p, q}^{m, n} (z) = \frac{1}{2 \pi i} \int\limits_{\mathcal{L}} \mathcal{H}_{p, q}^{m, n} z^s d s	
		\end{equation*}
		where $z^s = exp[s (ln |z| + i arg z)]$,
		\begin{equation*}
			\mathcal{H}_{p, q}^{m, n} = \frac{A (s) B (s)}{C (s) D (s)}, \qquad	
			A (s) = \prod_{j = 1}^m \Gamma(b_j - \beta_j s),\quad B (s) = \prod_{j = 1}^n \Gamma (1 - a_j + \alpha_j s),	
		\end{equation*}
		\begin{equation*}
			C (s)  = \prod_{j = m + 1}^q \Gamma(1 - b_j + \beta_j s),\quad  D (s) = \prod_{j = n + 1}^p \Gamma (a_j - \alpha_j s).	
		\end{equation*}
		Here $m,n,p,q$ are integers satisfying 
		$0 \le n \le p, \quad 1 \le m \le q, m^2 + n^2 \ne 0$,
		$a_j (j = 1, \dots, p), b_j (j = 1, \dots, q)$ are complex numbers.
		
		The integration contour $\mathcal{L}$ could be chosen in different ways:
		\begin{itemize}
			\item {$\mathcal{L} = \mathcal{L}_{- i \infty, i \infty}$ chosen to go from $- i \infty \quad$ to $\quad i \infty$} leaving to the right all poles of $\mathcal{P} (A)$  of the functions $\Gamma$ in $A (s)$ and to the left all poles of $\mathcal{P} (B)$  of the functions $\Gamma$ in $B (s)$;
			\item {$\mathcal{L} = \mathcal{L}_{i \infty}$ is a loop beginning and ending at $+ \infty$ and encircling in the negative direction all the poles of $\mathcal{P} (A)$;}
			\item {$\mathcal{L} = \mathcal{L}_{- i \infty}$ is a loop beginning and ending at $- \infty$ and encircling in the negative direction all the poles of $\mathcal{P} B)$.}
		\end{itemize}

		\subsubsection{Wright Functions}
		The Write function  is defined by the series representation that is convergent in the whole $z$-complex plane \cite{wf1,gor99,gor0,kil02}
		\begin{equation*}
			W_{\lambda, \mu} (z) = \sum_{n = 0}^{\infty} \frac{z^n}{n! \Gamma (\lambda n + \mu)}, \quad \lambda > - 1, \mu \in \mathcal{C}.
		\end{equation*}
		
		The integral representation of the Write function is written as 
		\begin{equation*}
			W_{\lambda, \mu} (z) = \frac{1}{2 \pi i} \int\limits_{Ha} e^{\sigma + z \sigma^{- \lambda}} \frac{d \sigma}{\sigma^{\mu}}
		\end{equation*}
		where $Ha$ is the Hankel path (a loop that starts from $- \infty$ along the lower side of the negative real axis, encircles the circular area the origin with radius $\epsilon \rightarrow 0$ in the positive sense, and ends at $- \infty$ along the upper side of the negative real axis).  
		
		There are Write-type {\em auxiliary} functions $	F_{\nu} (z) = W_{- \nu,0} (z)$, $ M_{\nu} (z) = W_{- \nu, 1 -\nu} (z)$,
		where $0 < \nu <1$; these functions are related $F_{\nu} (z) = \nu z M_{\nu} (z)$.  
		
		The series representations of the {\em auxiliary} functions are
		\begin{equation*}
			F_{\nu} (z) = \sum_{n = 1}^{\infty} \frac{(- z)^n}{n! \Gamma (- \nu n)} = \frac{1}{\pi} \sum_{n = 1}^{\infty} \frac{(- z)^{n - 1}}{n!} \Gamma (\nu n + 1) \sin (\pi \nu n)
		\end{equation*}
		and 
		\begin{equation*}
			M_{\nu} (z) = F_{\nu} (z) = \sum_{n = 0}^{\infty} \frac{(- z)^n}{n! \Gamma [- \nu n + (1 - \nu))]} = \frac{1}{\pi} \sum_{n = 1}^{\infty} \frac{(- z)^{n - 1}}{(n - 1)!} \Gamma (\nu n) \sin (\pi \nu n).
		\end{equation*}

		\section{Solution  of Fractional Differential Equations}
		\subsection{Analytical Methods}
		Numerous approximate analytical methods  are known:
		\begin{itemize}
			\item {the Adomian decomposition method (ADM) \cite{zha1};}
			\item {the combined Laplace-Adomian method (CLAM) \cite{waz10};}
			\item {the variational iteration method (VIM) \cite{y11,wu10,far11,yan13,yan13b}\footnote{
					VIM includes three steps to determine the variational iteration formula:
					\begin{enumerate}
						\item {establishing the correction functional;}
						\item {identifying the Lagrange multipliers;}
						\item {determining the initial iteration.}
					\end{enumerate} 	
					
					The second step is the crucial one \cite{zha9}.
					
				} and its local (LVIM) \cite{yang12a,yan13b,zha9,yang15} and fractional (using the fractional order Lagrange multipliers) \cite{khan,zha15b}  variants; }
			\item {the homotopy perturbation method (HPM)\cite{abb06,zha6,yil10,sin11a,raf12,wei17} and its modification \cite{yan13c} and local fractional variant (LFHPM) \cite{yang15a};}
			\item {the differential transformation method \cite{jon09,all09,gha15a};}
			\item {the heat-balance integral method (HBIM)\cite{chr10,hr11,chr11};}
			\item {the fractional complex transform method (FCTM) \cite{yang,y16,y16a,he12a,su13};}
			\item {the local fractional Fourier series method (FSM) \cite{yang13,zha14a}; }
			\item {the modified simple equation method \cite{jaw10,you14};}
			\item {the method of images (limited to special spatial symmetries);}
			\item {the Mellin integral transform method \cite{luc13};}
			\item {the local fractional decomposition method (LFDM) \cite{ahm15};}
			\item {the fractional sub-equation method \cite{y19,guo12};} 
			\item {the Sumudu transform methods \cite{wat02a,dem15} and its variant --- the local fractional homotopy perturbation Sumudu transform mehod \cite{zhao17};}
			\item {the theta-method \cite{asl14};}
			\item {the Picard succesive approximation method (PSAM) \cite{yang12b,yan14};}
			\item {the local Laplace transforms.}
		\end{itemize}
		
		Frequently analytical methods are variants of perturbation methods \cite{vd}).
		For example, He \cite{he1,he2,he3} based his  method to solve the general equation 
		$A (u) - f (r) = 0 $ 
		with the general differential operator $A$ divided into linear $L$ and nonlinear $N$ parts 
		$L (u) + N (u) - f (r) =0$ 	on the approach of Liao \cite{liao} (who used the two-parameter family of  equations) by considering the one-parameter family
		$(1 - p) L (u) +  p N (u) = 0$.
		
		He constructed the homotopy $v(r, p): \Omega \times [0, 1] \rightarrow R $ that satisfies
		\begin{equation*}
			H(v, p) = (1 - p) [L (v) - L (v_0)] + p [A (v) - f (r)] = 0
		\end{equation*}
		where the homotopy parameter $p \in [0, 1]$, 
		$v_0$ is the initial approximation.
		
		Evidently, $H (v, 0) = L (v) - L (v_0) =0$ 	and 
		$H (v, 1) = [A (v) - f (r)] = 0$. 
		
		In topology, $L (v) - L (v_0)$ is called deformation.
		The  homotopy parameter $p$ is considered as a small parameter and the solution is written as a series
		\begin{equation*}
			v = v_0 + p v_1 + p^2 v_2 + p^3 v_3 + \dots
		\end{equation*}
		and when $p \rightarrow 1$
		\begin{equation*}
			u = \lim_{p \rightarrow 1}v = v_0 + v_1 + v_2 + v_3 + \dots.
		\end{equation*}
		
		The Adomian decomposition method (ADM) \cite{ad1,ad2} does not use linearization, perturbation or the Green's functions. The accuracy of the approximate analytical solutions can be verified by the direct substitution.   
		
		The initial value is written as $L u + R u + N u = g$ 	where $L$ is the linear operator to be inverted, $R$ is the linear remainder operator and $N$ is the nonlinear operator.
		Thus $L^{- 1} L u = u - \Phi$, $\Phi$ incorporates the initial values.
		
		The solution and nonlinear term are decomposed into series
		\begin{equation*}
			u = \sum_{n = 0}^{\infty} u_n, \qquad
			N u = \sum_{n = 0}^{\infty} A_n
		\end{equation*}
		where $A_n$ are the Adomian polynomials for $N u = f (u)$ are \cite{ad2}
		\begin{equation*}
			A_n = \frac{1}{n!} \frac{\partial^n}{\partial \lambda^n} f \left(\sum_{k = 0}^{\infty} u_k \lambda^k  \right), \qquad n = 0, 1, 2, \dots.
		\end{equation*}
		
		Finally 
		\begin{equation*}
			\sum_{n = 0}^{\infty} u_n	= \Phi + L^{- 1} g - L^{- 1} \left[R  \sum_{n = 0}^{\infty} u_n + \sum_{n = 0}^{\infty} A_n \right].
		\end{equation*}
		
		The nonlinear term $Nu (x, t)$ can be also decomposed 
		\cite{yan13c} as
		\begin{equation*}
			N u = \sum_{n = 0}^{\infty} p^n H_n	 
		\end{equation*}
		where  He's polynomials are \cite{noor,gho09}
		\begin{equation*}
			H_n (u_0, u_1, \dots, u_n)	= \frac{1}{n!} \frac{\partial^n}{\partial p^n} \left[N \left( \sum_{i= 0}^n p^i u_i\right)  \right].
		\end{equation*}
		
		The fractional sub-equation method includes several steps 
		\cite{guo12}:
		\begin{itemize}
			\item {Transformation of the nonlinear fractional equation in two variables $x$ and $t$ $\quad  D_t^{\alpha} u,  D_x^{\alpha} u, \dots) = 0, \quad 0 < \alpha \le 1$,
			$D_t^{\alpha} u$ and $D_x^{\alpha} u$ are Jumarie modification of the RL derivtives,
			using the travelling wave transformation 
			$u (x, t) = u (\xi), \xi = x + c t$,
			where $c$ is a constant to be determined, to the equation 				\begin{equation}
					\label{g9}
					P (u, c u^{\prime}, u^{\prime}, c D_{\xi}^{\alpha} u,  D_{\xi}^{\alpha} u, \dots) = 0.
				\end{equation}	
			}
			\item {The solution of the equation (\ref{g9}) is assumed to have the form
				\begin{equation*}
					u (\xi) = \sum_{i = - n}^{- 1} a_i \phi^i + a_0 + \sum_{i = 1}^{n} a_i \phi^i,
				\end{equation*}
				where $a_i (i =  - n, - n + 1, \dots, n - 1, n)$ are constants to be determined, $\phi = \phi (\xi)$ are functions that satisfy the following Riccati equation $D_{\xi}^{\alpha} \phi (\xi) = \sigma \phi^2 (\xi)$,
				$\sigma$ is a constant.
			} 
			\item {Formulation of a set of overdetermined nonlinear algebraic equations for $c$ and $a_i (i =  - n, - n + 1, \dots, n - 1, n)$ \cite{guo12}. }
		\end{itemize}

		\subsection{Numerical Methods}
		Diethelm et al. \cite{die06} listed the requirement to the numerical methods that should be 
		{convergent,} {consistent of some reasonable order $h^p$,} {stable,} {reasonably inexpensive to run,} {reasonably easy to program.}
			
		Numerous methods are used in practice: finite difference, finite elements, radial basis functions, spectral methods, meshfree  methods.
		The numerical methods for the fractional differential equations usually are constructed by the modification of the methods for the ordinary differential equations but require significantly more computation time and storage. The approximation of the fractional derivative needs the computation of the convolution integral that requires to sample and multiply the behaviour of two functions over the whole of the interval of integration leading to the operation count of $O (n^2)$ where $n$ is number of sampling points \cite{ford}.
		
		The reduction of the computational efforts is related to the fading memory property of the fractional derivatives that allows to restrict the integration interval --- using the {\em short memory principle} \cite{den07a,den15} (also {\em fixed memory principle}  \cite{ford} and {\em logarithmic memory principle} \cite{han11}), and using adaptive time stepping and basis selection \cite{bru10}.
		
		Numerous methods are used to solve the fractional differential equations in practice: the finite difference \cite{fd1,fd2} (both the explicit, e.g. Euler \cite{odi08} and the implicit \cite{impl1,impl2}, e.g., the Crank-Nicolson \cite{cn1,cn2} or the  alternating direction implicit \cite{zha16a,zha19b}  schemes, compact schemes \cite{cfd1,cfd2,cfd3,cfd4}), the finite elements \cite{hua08,jia11,zeng15,fem1,fem2} (including least squres FEM \cite{fix}, Galerkin FEM \cite{sin06b,jin15}, discontinuous Galerkin FEM \cite{mus13}),  the spectral methods \cite{can06,esm11,doha11,zay14}, the meshfree  methods \cite{dua08,shi12} (including the radial basis functions methods that exploit cubic $\phi = r^3$, Gaussian $\phi = exp (- r^2 / c^2)$,  multiquadrics $\phi = \sqrt{c^2 + r^2}$  or inverse  multiquadrics $\phi =1 / \sqrt{c^2 + r^2}$ functions \cite{rb1,rb2,rb3}), Legendre wavelet collocation method \cite{hey14}.
		
		Bahuguna et al. \cite{bah09}, Hanert \cite{ha10}, Deng et al. \cite{den04} and Deng \& Li \cite{deng12}, Ford \& Simpson \cite{ford}, Ford \& Connolly \cite{ford1}, Momani et al. \cite{mom} reported the results of the comparison of several numerical methods.

\end{document}